\documentclass[preprint]{elsarticle}
 \biboptions{longnamesfirst,angle,semicolon}
\usepackage{amsmath,amsthm,amsfonts,amssymb,tensor}

\usepackage{atbegshi}
\AtBeginDocument{\AtBeginShipoutNext{\AtBeginShipoutDiscard}}

\makeatletter
\DeclareRobustCommand\frownotimes{\mathbin{\mathpalette\frown@otimes\relax}}
\newcommand{\frown@otimes}[2]{%
  \vbox{
    \ialign{##\cr
      \hidewidth$\m@th#1{}_^$\kern-\scriptspace\hidewidth\cr
      \noalign{\nointerlineskip\kern-1pt}
      $\m@th#1\otimes$\cr
    }%
  }%
}
\makeatother

    \makeatletter
    \def\ps@pprintTitle{%
       \let\@oddhead\@empty
       \let\@evenhead\@empty
       \let\@oddfoot\@empty
       \let\@evenfoot\@oddfoot
    }
    \makeatother

\usepackage{multicol}
\usepackage{color}
\usepackage{xspace}
\usepackage{enumerate}
\setcitestyle{square}

\DeclareRobustCommand{\LaTeX}{L\kern-.26em%
        {\sbox\z@ T%
         \vbox to\ht\z@{\hbox{\check@mathfonts
           \fontsize\sf@size\z@
           \math@fontsfalse\selectfontd
          A\,}%
         \vss}%
        }
     \kern-.15em%
    \TeX}
\makeatother

\newtheorem{thm}{Theorem}[section]

  \newtheorem{remark}[thm]{Remark}
  \newtheorem{corollary}[thm]{Corollary}
 \newproof{pf}{Proof}
 \numberwithin{equation}{section}
\begin{document}
\title{On an isomorphism theorem for the Feichtinger's Segal algebra on locally compact groups\tnoteref{t1}}
 \author[]{R. Lakshmi Lavanya\tnoteref{t2}}
 \tnotetext[t1] {Dedicated to the memory of Prof. Elias M. Stein}
  \tnotetext[t2]{The author's research is supported by DST-SERB, Government of India, under the Early Career Research Award Scheme(No. ECR/2017/003366).}

%


 \address{Indian Institute of Science Education and Research (IISER) Tirupati,\\
 Karakambadi Road, Mangalam,\\
 Tirupati - 517507, India.}
 \begin{abstract}
  In this article we observe that a locally compact group $G$ is completely determined by the algebraic properties of its Feichtinger's Segal algebra $S_0(G).$ Let $G$ and $H$ be locally compact groups. Then any linear (not necessarily continuous) bijection of $S_0(G)$ onto $S_0(H)$ which preserves the convolution and pointwise products is essentially a composition with a homeomorphic isomorphism of $H$ onto $G.$
 \end{abstract}
\vspace*{-.5pc}
 
\maketitle

\section{Introduction} 
\setcounter{equation}{0}
\label{Prelims}

In 1981, Feichtinger introduced and studied a new Segal algebra on locally compact Abelian groups(\cite{Fei}). This algebra, now called the Feichtinger algebra has been extensively studied over the years. A recent survey of this progress can be found in \cite{Jak}. The genesis of the study of the Feichtinger's Segal algebra for locally compact groups came up almost a quarter century later(\cite{Sp}). In \cite{Sp}, the author introduces an operator space structure on Feichtinger's Segal algebra for locally compact groups. It was also observed that $S_0(G)$ for a locally compact group has many functorial properties very similar to that in the Abelian setting. Infact, some of these properties could be obtained even without exploiting the operator space structure of the Feichtinger's Segal algebra. One such example is the isomorphism theorem, viz., Theorem \ref{Spthm} proved in \cite{Sp}. This result states that any bounded linear bijection which preserves the convolution and pointwise products on the Feichtinger's Segal algebra determines the underlying locally compact group.\\

The aim of this article is to observe that the algebraic properties of the Feichtinger's Segal algebra determine the underlying locally compact group. More precisely, a locally compact group is completely determined by any linear bijection (not necessarily continuous), which preserves the convolution and pointwise products on its Feichtinger's Segal algebra.\\

For more on isomorphism theorems, we refer to \cite{KA}, \cite{Sp}, \cite{Wa}, \cite{We} and the references therein.

\section{Preliminaries and Notation}

For a locally compact group $G,$ let $\mathcal{C}_c(G)$ denote the space of all compactly supported continuous functions on $G.$ As per the standard convention, we denote by $A(G)$ and $B(G),$ the Fourier and the Fourier-Stieltjes algebras on $G,$ respectively. The Feichtinger's Segal algebra on $G,$ denoted by $S_0(G),$ is defined as follows(\cite{Sp}):\\

Let $\mathcal{A}$ be a regular Banach algebra with Gelfand spectrum $X.$ An ideal $\mathcal{I}$ of $\mathcal{A}$ is said to have \textit{compact support} if $ supp \ \mathcal{I}  =\overline{X\setminus h(\mathcal{I})}$ is compact in $X.$ Here, $h(\mathcal{I})$ denotes the \textit{hull} of $I,$ defined as $h( \mathcal{I} ) = \{x\in X: u(x)=0 \textrm{ \  for  \ all \ } u \in  \mathcal{I}  \}.$

Let $ \mathcal{I} $ be a non-zero ideal in $A(G)$ which has compact support. Identify $l^1(G)$ with the closed linear span of the Dirac measures $\{\delta_s: s\in G\}$. Consider $ l^1(G)$ with the maximal operator space structure. 

\noindent Define  $$q_ \mathcal{I} : l^1(G){\widehat{\otimes}} \mathcal{I}  \rightarrow A(G)$$as
$$q_ \mathcal{I} (x\otimes u)  = \sum \limits_{s\in G} \alpha_s\ s \cdot u,$$
where $x= \sum\limits_{s\in G} \alpha_s \delta_s$ and $s\cdot u(t) = u(s^{-1}t), \ t\in G.$
Here, $\widehat{\otimes}$ denotes the operator projective tensor product. The \textit{Feichtinger's Segal algebra} on $G$ is defined as
$$S_0(G) = ran(q_ \mathcal{I} ) \cong l^1(G) \widehat{\otimes} I \Big/ ker(q_ \mathcal{I} ).$$
Equipped with the quotient operator space structure, $S_0(G)$ becomes an operator space. \\

For functions $u,v \in S_0(G), \ u\ast v $ denotes the usual convolution product, given by 
$$u \ast v(x) = \int\limits_G u(y) \ v(y^{-1}x) \ dy.$$
Here, we use $dy$ to denote a fixed left Haar measure on $G.$\\

In \cite{Sp}, it was illustrated that $S_0(G)$ satisfies several functorial properties analogous to its Abelian counterpart.\\

The following isomorphism theorem for $S_0$ was proved in \cite{Sp}.

\begin{thm}[\cite{Sp}]\label{Spthm}
Let $G$ and $H$ be locally compact groups and $ \Phi: S_0(G) \rightarrow ~S_0(H)$ be a bounded linear bijection which satisfies 
$$\Phi(uv) = \Phi u \ \Phi v \textrm{ \ and \ } \Phi(u\ast v) = \Phi u \ast \Phi v$$
for every $u,v\in S_0(G).$ There is a homeomorphic isomorphism $\alpha:H \rightarrow G$ such that $$\Phi u = u\circ \alpha$$
for each $u$ in $S_0(G).$
\end{thm}

The aim of the present article is to observe that the conclusion of the above result holds good even if the map $\Phi$ is not apriori continuous (see Theorem \ref{LASG}). This leads us to the fact that the algebraic properties of the Feichtinger's Segal algebra solely determine the underlying locally compact group.\\
 
\noindent  For further details and properties of $S_0(G),$ we refer to \cite{Sp}.


\section{An Isomorphism theorem for $S_0(G)$}
First we obtain the form of any bijection which is multiplicative on $S_0(G).$ Later we will obtain the isomorphism theorem by additionally assuming that this bijection preserves the convolution product too.\\

In this section, we use letters $x,y, \cdots$ to denote the elements of the respective groups, as indicated in the context. 

Our preliminary result is the following, which exploits only the pointwise product structure of $S_0(G).$
\begin{thm}\label{Lamult}
Let $G $ and $H$ be locally compact groups and $\Phi:S_0(G) \rightarrow S_0(H)$ be a  bijection which satisfies the following conditions for every $u,v \in S_0(G):$ 
\begin{description}
\item[(1)] $\Phi(u+\overline{v}) = \Phi (u) + \overline{\Phi(v)}$
\item[(2)]$ \Phi (uv) = \Phi u \ \Phi v$
\end{description} 
Then there exists a decomposition of $H$ as $H=H_1 \cup H_2$ and a homeomorphism $\alpha:H\rightarrow G$ such that 
$$\Phi u = \chi_{H_1} \ (u \circ \alpha) + \chi_{H_2}\  (\overline{u\circ \alpha}),$$
for each $u \in S_0(G).$
\end{thm}

\begin{pf} The proof consists of 7 steps. Some of these overlap with a characterisation of the Fourier tranform on the Fourier twins on locally compact Abelian groups (as in Theorem 3.2 of \cite{Lamona}). For the sake of completeness, a proof of these ideas (steps 1-5) is included in the Appendix. Here, we only state these steps.\\

\noindent \textbf{Step 1.} If $u\in \mathcal{C}_c(G),$ then $\Phi u \in \mathcal{C}_c(H).$ \\

\noindent \textbf{Step 2.} For any $x_0 \in G,$ there exists $y_0 \in H$ such that $\Phi u (y_0) = 0$ whenever $u(x_0) = 0,$ for any $u \in S_0(G).$\\
\noindent We say that a function $u:G\rightarrow \mathbb{C}$ satisfies condition $(\star)$ if the following holds good:
$$(\star) \ \ \ \ \ \ \ \ \ \ u(x_0)=0 \textrm{ \ if \ and\  only \ if \ } \Phi u(y_0)=0.\hspace{4cm}$$
\noindent \textbf{Step 3.} Define a map $\beta : G \rightarrow H$ as follows: $\beta(x) = y$ if $\Phi u (y) 
= 0$ whenever
$u(x) 
= 0,$ for any $u\in S_0(G).$ Then the map $\beta$ is well-defined.\\

\noindent \textbf{Step 4.} The map $\beta : G \rightarrow H$ is a bijection.\\

\noindent \textbf{Step 5.} The map $\beta$ is a homeomorphism of $G$ onto $H.$\\

\noindent \textbf{Step 6.} 'Extension' of the map $\Phi $ to scalars.\\

\noindent \textit{Illustration of Step 6.} For $u,v \in
S_0(G),$ and $c (\neq 0)\in \mathbb{C},$ we have
$$ \Phi (c u)(x)\ \Phi v(x) = \Phi (c uv)(x) = \Phi (u)(x) \ \Phi (c v)(x) , \ x\in G.$$
Let $w \in S_0(G)$ be such that $\Phi w(x) \neq 0$ for any $x\in
G.$ Then we have \noindent \begin{eqnarray*}
 \Phi (c u)(x)  &=& \frac{\Phi (c w)(x)}{\Phi w(x)} \ {\Phi u}(x) \textrm{ \ for \ all \ } u\in S_0(G)\\
   &=& m(c ,x) \ \Phi u(x) \textrm{ \ (say)}.
\end{eqnarray*}

\noindent When $c=0,$ we have $m(c,x) \Phi u(x) = \Phi(cu)(x)=0,$ for all $x\in H$ and for all $u \in S_0(G).$ Hence it is reasonable to define $m(0,x) = 0$ for all $x\in H.$\\
\noindent (Here, we have used that $\Phi(0_G)=0_H,$ where $0_G$ and $0_H$ are the zero functions on the respective groups. This follows by condition(1) since $\Phi(0_G)=\Phi(0_G+0_G)=\Phi(0_G)+\Phi(0_G).$)
$$ \textrm{Thus \hspace{3cm} }\Phi (c u)(x) = m(c ,x)\ \Phi u(x), \textrm{ \  for\  all \ }x\in H. \hspace{4cm}$$ In particular, the function $m(\cdot,\cdot)$ is continuous in the
second variable, as a function
of $x\in H.$\\

\noindent \textbf{Step 7.}  For fixed $y\in H,$ the map $m(\cdot,y) :\mathbb{C} \rightarrow
\mathbb{C}$ is an additive and multiplicative bijection, which
maps $\mathbb{R}$ onto $\mathbb{R},$ and hence we have either
$m(c,y) = c,$ for all $c \in \mathbb{C},$ or $m(c,y)
= \overline{c},$ for all
$c \in \mathbb{C}.$\\

\noindent \textit{Proof of Step 7.} 
Choose $u\in S_0(G),$ with $u(x) \neq 0$ for any $x\in G.$
Then by the condition
$(\star),$ we get $\Phi u(\beta(x))  \neq 0,$ for any $x\in G.$\\

Suppose $m(c,y) =m(d,y)$ for some $c,d \in
\mathbb{C}.$ Then
$$\Phi(c u) (y) = m(c,y) \ \Phi u(y)= m(d,y) \ \Phi u(y) =
\Phi(d u)(y) , \ y\in H.$$ As $\Phi$ is a bijection, this gives $cu=du.$ Since $u$ is never vanishing, we get that $c=d.$\\

\noindent By hypothesis(1), we have \noindent \begin{eqnarray*}
 m(c +\overline{d},y)\ \Phi u(y)  &=& \Phi((c +\overline{d})u)(y)
=\Phi(c  u +\overline{d}u)(y) \\
   &=& \Phi(c u)(y) + \overline{\Phi(d \overline{u})}(y) = (m(c,y ) 
   +\overline{m(d,y )}) \ \Phi u(y).
\end{eqnarray*}
Since $\Phi u$ is never zero, we get $m(c +\overline{d},y) =
m(c,y )+\overline{m(d,y)}.$ In particular,
$m(d,y) \in \mathbb{R}$ whenever $d \in \mathbb{R}.$\\

\noindent Now, hypothesis(2) gives
$$m(c d , y) \Phi u(y) = \Phi(cd u)(y) = 
m(c,y)\Phi(d u)(y) = m(c,y ) \ m(d,y) \Phi u(y).$$
Again, since $\Phi u$ is nowhere vanishing, we get $m(c
d,y)=m(c ,y) \ m(d,y)$
for all $c,d \in \mathbb{C}.$

Thus $m(\cdot,y)$ is a field automorphism of $\mathbb{C}$ which preserves $\mathbb{R}.$ So
for fixed $y\in H,$ we have either $m(c,y) = c$ for all $c\in \mathbb{C},$ or $m(c,y) = \overline{c}$ for all $c\in \mathbb{C}.$ Let
\begin{eqnarray*}
H_1&=&\{y\in H: m(c,y) = c \textrm{ \ for \ all \ } c\in \mathbb{C}\},\\
H_2&=&\{y\in H: m(c,y) = \overline{c} \textrm{ \ for \ all \ } c\in \mathbb{C}\}.
\end{eqnarray*}

\noindent \textbf{Step 8.} For $u\in S_0(G),$ and $x_0 \in
G,$ we have
$\Phi u[\beta (x_0)] = m(u(x_0),\beta(x_0)).$\\

\noindent \textit{Proof of Step 8.} As before, choose $u \in
S_0(G)$ such that $u(x) \neq 0$
for any $x \in G.$ Then $\Phi u(y) \neq 0$ for any $y \in H.$ Let $v\in S_0(G).$ \\

\noindent Define $$w(x): = v(x_0) \ u(x) - v(x) \ u(x), \ x\in
G.$$ Then $w\in S_0(G)$ and $w(x_0) = 0.$ By Condition
$(\star),$ we have $\Phi w[\beta (x_0)] =0.$ This gives \noindent
\begin{eqnarray*}
0=\Phi w [\beta (x_0)] &=& \Phi[v(x_0) \cdot u - v \cdot u][\beta (x_0)]\\
 &=& m(v(x_0),\beta(x_0)) \
\Phi u[\beta (x_0)] - \Phi v[\beta (x_0)] \ \Phi u[\beta (x_0)].
\end{eqnarray*}
 Since $\Phi u$ is never zero, this gives
$\Phi v[\beta (x_0)] = m(v(x_0),\beta(x_0)).$\\

\noindent Since $\alpha =\beta ^{-1},$ using Step 7, we get that
$$\textrm{either \ } \Phi v(x_0) = v[\alpha (x_0)] \textrm{ \ or \ } \Phi v(x_0) = \overline{v[\alpha (x_0)]}.$$

\noindent Combining this with the definition of the sets $H_1$ and $H_2$ gives the representation $\Phi v =\chi_{H_1} \  \alpha \circ v+\chi_{H_2} \ \overline{\alpha \circ v}. $ \qed\\
\end{pf}


\begin{corollary}
In addition to the hypothesis of Theorem \ref{Lamult}, if H is connected, then the following conclusion holds:
\\There exists a homeomorphism $\alpha:H\rightarrow G$ such that 
$$\textrm{either \ } \Phi u =u \circ \alpha \textrm{ \ for \ all \ } u \in S_0(G), $$
$$ \textrm{or \ } \Phi u =  \overline{u\circ \alpha}\textrm{ \ for \ all \ } u \in S_0(G). $$
\end{corollary}
\begin{pf}
Since the map $m(\cdot,\cdot): \mathbb{C} \times H \rightarrow \mathbb{C}$ is continuous in the second variable, we get by Step 7 of the previous theorem that either $m(c,y) = c$ for all $y\in H, $ or $m(c,y) = \overline{c}$ for all $ y\in H.$\qed\\
\end{pf}

As of now, we have not made use of the convolution product in $S_0(G).$ As mentioned in \cite{Sp}, it is reasonable to expect that any isomorphism theorem for the Feichtinger's Segal algebra involves both the  convolution and the pointwise product. This leads us to the following result. The proof of this result is along the lines of the proof of Step 6 of Theorem 3.2 in \cite{Lamona}. We state it here for the sake of completeness.

\begin{thm}\label{Laconv}
In addition to the hypothesis of Theorem \ref{Lamult}, assume that $\Phi$ satisfies 
$$\Phi(u\ast v) = \Phi u\ast \Phi v, \textrm{  \ for \ all \  }u,v \in S_0(G).$$ Then there exists a decomposition of $H$ as $H=H_1 \cup H_2$ and a homeomorphic isomorphism $\alpha:H\rightarrow G$ such that 
$$\Phi u = \chi_{H_1} \ (u \circ \alpha) + \chi_{H_2}\  (\overline{u\circ \alpha}),$$
for each $u \in S_0(G).$
\end{thm}

\begin{pf} In view of the proof of Theorem \ref{Lamult}, we are left with proving that $\alpha$ is a homomorphism.\\

\noindent \textit{Claim:} The map $(\alpha^{-1}=)\beta: G \rightarrow H$  is a homomorphism.\\

\noindent \textit{Proof of Claim.} Suppose $\beta (xy) \neq \beta (x) \beta (y)$ 
for some $x,y\in G.$\\

\noindent Then there exist disjoint neighbourhoods
$V_{xy},V_{x\diamond y}$ in $H$ with $\beta (xy) \in ~V_{xy}$ and
$\beta (x) \beta (y) \in V_{x\diamond y }.$ By continuity of
the map $\beta ,$ this gives rise to a neighbourhood $W_{xy,1}$
of $xy$ in $G$ with $\beta (W_{xy,1}) \subseteq V_{xy}.$ By continuity
of multiplication in $G,$ we get neighbourhoods $W_{x,1}, W_{y,1}$
of $x$ and $y,$ respectively, in $G,$ such that $W_{x,1}W_{y,1} \subseteq
W_{xy,1}.$ Thus \noindent \begin{eqnarray} \vspace{-1cm} \beta
(W_{x,1} W_{y,1}) &\subseteq& \beta (W_{xy,1}) \subseteq
V_{xy}.\label{eq1}
\end{eqnarray}

On the other hand, by continuity of multiplication in $H,$
$\beta( x) \beta (y) \in V_{x\diamond y }$ gives
neighbourhoods $V_{x,2},V_{y,2}$  in $H,$ such that \noindent
\begin{eqnarray} \beta (x)\in V_{x,2}, \ \beta (y) \in V_{y,2}
\textrm{ \ and \ } V_{x,2}\ V_{y,2}\subseteq V_{x\diamond y
}.\label{eq2}
\end{eqnarray}
 This
implies there exist neighbourhoods $W_{x,2}, W_{y,2}$ of $x$ and
$y,$ respectively, with $\beta (W_{x,2})\subseteq V_{x,2}$ and
$\beta (W_{y,2}) \subseteq V_{y,2}.$\\

Define $W_x= W_{x,1}\cap W_{x,2}, \ W_y = W_{y,1}\cap W_{y,2}.$
Then \noindent \begin{eqnarray}
  \beta (W_x)& \subseteq & \beta (W_{x,2}) \subseteq V_{x,2} =V_x \ (say) \nonumber \\
   \beta (W_y)& \subseteq & \beta (W_{y,2}) \subseteq V_{y,2} =V_y \ (say)\nonumber \\
   \beta (W_x W_y)& \subseteq & \beta (W_{x,1} W_{y,1})  \subseteq \beta (W_{xy,1})
    \subseteq V_{xy} \label{eq3}
\end{eqnarray}
Choose $u_x,u_y \in S_0(G)$ such that $Supp \ u_x \subseteq
W_x, Supp\ u_y \subseteq W_y$ and $u_x \ast ~u_y \not \equiv~0.$
Let $v_x = \Phi u_x$ and $v_y = \Phi u_y.$ Then $\Phi (u_x\ast u_y) = v_x \ast
v_y \not\equiv 0.$\\

\noindent We have \noindent \begin{eqnarray}
 \nonumber Supp(v_x \ast v_y) &=&Supp \ \Phi (u_x \ast u_y)  \subseteq  (Supp \ \Phi u_x) \ (Supp \ \Phi u_y) \\
\nonumber &= & \beta (Supp \ u_x) \ \beta (Supp \ u_y)
\subseteq \beta( W_x) \ \beta (W_y)
\\ & \subseteq & V_{x,2}\ V_{y,2} \subseteq V_{x\diamond y } \textrm{   \ \ (by \
(\ref{eq2}))} \label{eq4}
    \end{eqnarray}
But $Supp \ (u_x  \ast u_y) \subseteq (Supp \ u_x) \ (Supp \ u_y)
\subseteq W_x\ W_y.$ By $(\ref{eq3}),$ this gives
 \noindent \begin{eqnarray}
 \nonumber Supp(v_x \ast v_y) &= &Supp (\Phi u_x \ast \Phi u_y) = Supp \ \Phi (u_x \ast u_y)\\
  & =& \beta (Supp (u_x  \ast u_y)) \subseteq  \beta (W_x \ W_y)\subseteq
 V_{xy}. \label{eq5}
 \end{eqnarray}
From $(\ref{eq4})$ and $(\ref{eq5})$, we get $$Supp(v_x \ast v_y)
 \subseteq V_{x\diamond y } \cap V_{xy} = \emptyset.$$ This gives $v_x\ast
v_y \equiv 0,$ a contradiction.
This proves the multiplicativity of the map $\beta .$\\

Thus the map $\beta $ is a multiplicative homeomorphism of G
onto $H.$ As $\alpha = \beta^{-1},$ this proves our result.\qed\\
\end{pf}

We now state an isomorphism theorem for the Feichtinger's Segal algebra.
\begin{thm}\label{LASG}
Let $G $ and $H$ be locally compact groups and $\Phi:S_0(G) \rightarrow ~S_0(H)$ be a  linear bijection which satisfies the following conditions for every $u,v \in ~S_0(G):$ 
$$ \Phi (uv) = \Phi u \ \Phi v \textrm{ \ and \ } \Phi(u \ast v) = \Phi u \ast \Phi v.$$
Then there exists an isomorphic homeomorphism $\alpha:G\rightarrow H$ such that 
$$\Phi u = u \circ \alpha, \textrm{ \ for \ all \ } u \in S_0(G).$$
\end{thm}
\begin{pf}
We proceed as in Theorem \ref{Lamult}, till Step 5. \\

Let $u \in
S_0(G)$ be such that $u(x) \neq 0$
for any $x \in G.$ Then $\Phi u(y) \neq 0$ for any $y \in H$ (by Condition ($\star$)). Let $v\in S_0(G).$
$$\textrm{Define} \hspace{2cm} w(x): = v(x_0) \ u(x) - v(x) \ u(x), \ x\in
G. \hspace{4cm}$$ Then $w\in S_0(G)$ and $w(x_0) = 0.$ By Condition
$(\star),$ we get $\Phi w[\beta (x_0)] =0.$ Since $\Phi$ is linear, this gives \noindent
\begin{eqnarray*}
0=\Phi w [\beta (x_0)] &=& \Phi[v(x_0) \cdot u - v \cdot u][\beta (x_0)]\\
 &=&v(x_0)\Phi u [\beta(x_0)] -\Phi v[\beta(x_0)]\ \Phi u [\beta(x_0)].
\end{eqnarray*}
 Since $\Phi u$ is never zero, this gives
$\Phi v[\beta (x_0)] = v(x_0),$ whenever $v(x_0)\neq 0.$\\

\noindent If $v(x_0)= 0,$ then condition ($\star$) gives that $\Phi v(\beta(x_0)) = 0.$

\noindent Since $\alpha =\beta ^{-1},$ this proves that
$$\Phi v= v \circ \alpha \textrm{ for \ all \ } v \in S_0(G).$$

As in the proof of Theorem \ref{Laconv}, we get that $\alpha$ is an isomorphism.
\qed\\
\end{pf}

\begin{remark}
As stated in \cite{Sp}, the above isomorphism theorems too, hold good for any regular Banach algebra of functions in G (whose spectrum apriori need not be G), which is contained in $C_0(G).$ 
\end{remark}
\begin{remark}
In particular, Theorem \ref{LASG} holds good for the Rajchman algebra $B_0(G) = B(G) \cap C_0(G),$ defined  and studied in \cite{KA}.
\end{remark}

%
%

\section{Appendix}

\noindent \textbf{Part of proof of Theorem \ref{Lamult}:} (as in \cite{Lamona})\\

\noindent For $x\in G,$ let $C(x) = \{f\in S_0(G): f(x) \neq 0\}.$\\

\noindent \textbf{Step 1.}  If $u \in \mathcal{C}_c(G),$ then $\Phi u \in \mathcal{C}_c(H).$\\

\noindent \textit{Proof of Step 1.} If $v=1$ on $Supp \ u,$ we
have $u\cdot v =u.$ This gives $\Phi u= \Phi(u\cdot v) = \Phi u \cdot \Phi v,$
and so $\Phi v=1$ on the set $\{x : \Phi u(x) \neq 0\}.$ By the continuity of $\Phi v,$ we get $\Phi v=1$ on $Supp(\Phi u).$\\

For $u(\not \equiv 0)\in
\mathcal{C}_c(G),$ choose $v\in
S_0(G)$ such that $v=1$ on $Supp \ u.$ Then $\Phi v=1$ on
$Supp \ \Phi u.$
Since $\Phi v\in S_0(H)\subset C_0(H),$ we get $Supp \ \Phi u$ is compact.\\
\vspace{0.75cm}

\noindent \textbf{Step 2.} For any $x_0 \in G,$ there exists $y_0
\in H$ such
that $\Phi u\in C(y_0)$ whenever $u \in C(x_0).$\\

\noindent \textit{Proof of Step 2. }Let $E:=\{u\in S_0(G):
u(x_0) \neq 0\}.$ Fix a function $v\in \mathcal{C}_c(G)$ with
$v(x_0) \neq 0.$ By Step 1, we have $K: = Supp \ \Phi v$ is compact.\\

\noindent  For $u\in E,$ define $K_u:= K \cap Supp \ \Phi u.$ For
finitely many functions $u_0:=v,u_1,\cdots ,u_k\in E,$ we have
$\prod_{j=0}^k \ u_j\not \equiv 0,$ and so $\prod_{j=0}^k \ \Phi u_j=
\Phi (\prod_{j=0}^k \ u_j) \not \equiv 0,$ which gives $\cap _{j=0} ^k
\ K_{u_j} \neq \emptyset.$ This means, the collection $\{K_u: u\in
E\}$ of closed subsets of $K$ has finite intersection property.
Since $K$ is compact, this gives $\bigcap \limits_{u\in E} \ K_u
\neq \emptyset.$ Let $y_0 \in \bigcap
\limits_{u\in E} \ K_u.$\\

\noindent \textbf{Claim.} For $u\in S_0(G),$ we have $u(x_0) = 0$ if and only if $\Phi u(y_0)=0.$\\

\noindent \textit{Proof of Claim.} We prove the claim in two
separate
cases.\\

\noindent Suppose $u(x_0) \neq 0.$\\
 Then
$u$ never vanishes on a neighbourhood, say, $V$ of $x_0.$ Let
$v\in S_0(G)$ be such that $u\cdot v = 1$ on $V.$ Choose
$w\in S_0(G)$ such that $w=1$ on a neighbourhood $W$ of
$x_0,$ and satisfies $W\subseteq Supp \ w \subseteq V.$ Since
$u\cdot v=1$ on $Supp \ w,$ by Step 1, we get $\Phi (u\cdot v)=\Phi u\cdot
\Phi v =1$ on $Supp\
\Phi w.$ Since $h\in E,$ by definition, $y_0 \in Supp \ \Phi w .$ 
This implies $\Phi u(y_0) \neq 0,$ and hence $\Phi u\in C(y_0).$\\

\noindent Note that all our arguments till now can be applied to
the map $\Phi ^{-1}$ as well, and so we have proved that $u(x_0)\neq
0$ if and
only if $\Phi u(y_0) \neq 0.$\\

\noindent A function $u:G\rightarrow \mathbb{C}$ is said to
satisfy the condition $(\star)$ if the following holds:
$$\hspace{3cm}u(x_0)=0 \textrm{ \ if \ and\  only \ if \ } \Phi u(y_0)=0. \ \hspace{2cm} \ \ \ (\star)$$
By the above discussion, we have that all functions in
$S_0(G)$ satisfy the Condition $(\star).$\\

\noindent \textbf{Step 3.} Define a map $\beta : G \rightarrow
H$ as follows: $\beta (x)=y$ if $\Phi u\in C(y)$
whenever $u\in C(x).$ Then the map $\beta $ is well-defined.\\

\noindent \textit{Proof of Step 3.} Suppose for some $x_0 \in G,$
we have $\beta (x_0) = y_1$ and $\beta (x_0) = y_2$ with
$y_1\neq y_2.$ Let $V_1$ and $V_2$ be disjoint open neighbourhoods in $H,$ of
$y_1$ and $y_2,$ respectively. There exists functions $v_1$ and
$v_2$ in $S_0(H)$ which are supported
 in $V_1$ and $V_2,$ respectively, such that $v_1(y_1) \neq 0$ and $v_2(y_2) \neq 0.$ 
 Let $u_1,u_2 \in S_0(G)$ with $\Phi u_1=v_1$ and
 $\Phi u_2=v_2.$ Then $0\equiv v_1\cdot v_2 =\Phi u_1 \cdot \Phi u_2 = \Phi (u_1 \cdot u_2)$ and
  so $u_1\cdot u_2 \equiv
 0.$\\

 On the other hand, as $v_j(y_j) = \Phi u_j(y_j)\neq 0$ for $j=1,2,$ we
 have by
Condition $(\star)$ that $u_j(x_0) \neq
 0$
 for $j=1,2,$ which is in contradiction to $u_1\cdot u_2\equiv 0.$\\
\vspace{0.75cm}

\noindent \textbf{Step 4.} The map $\beta :G \rightarrow H$ is a bijection.\\

\noindent \textit{Proof of Step 4.} The hypotheses of the theorem
hold good for the map $\Phi ^{-1}$ as well. Applying the preceding
steps to the map $\Phi ^{-1}$ gives rise to a well-defined function,
say,
$\alpha :H \rightarrow G.$ Then $\alpha  = \beta ^{-1}$, proving that $\beta $ is bijective.\\
\vspace{0.75cm}

\noindent \textbf{Step 5.} The map $\beta $ is a homeomorphism of $G$ onto $H.$\\

\noindent \textit{Proof of Step 5.} First we prove that $\beta$
is continuous. Suppose not. Then there exists $x\in G,$ and a net
$(x_\tau)$ in $G$ with $x_\tau
\rightarrow x,$ but $\beta( x_\tau)$ does not converge to $\beta (x).$\\

\noindent Let $V$ be a neighbourhood of $\beta (x)$ in $H,$ such that
$\beta( x_\tau) \not\in V$ for any $\tau.$ Let $w\in
S_0(H)$ with $Supp \ w\subseteq V,$ and $w[\beta (x)] = 1.$
Let $v\in S_0(G)$ be such that $\Phi v=w.$ Then $\beta (x_\tau
)\not \in Supp \ \Phi v$ for any $\tau,$ and so $x_\tau \not\in Supp \
v$ for any $\tau.$ This gives
 $v(x_\tau) = 0$ for all $\tau,$ implying $v(x)= 0,$ which is not
possible by Condition $(\star),$ since $\Phi v[\beta (x)] =1.$\\

We observe that the above argument holds good when the maps $\Phi $
and $\beta $ are replaced with $\Phi ^{-1}$ and $\beta ^{-1}$
respectively, yielding that the map
$\beta : G \rightarrow H$ is a homeomorphism.\\

\noindent Since $\beta$ is a homeomorphism, Condition$(\star)$ gives:
$$Supp \ \Phi u = \beta (Supp \ u),\textrm{ \ for \ all \ } u\in S_0(G).$$
\begin{center}
---
\end{center}

\section*{References}
\bibliographystyle{spmpsci}      

\begin{thebibliography}{}
%


\bibitem{ER} E.G. Effros, Z.-J. Ruan, \textit{Operator Spaces,} London Math. Soc. Monogr. (N.S.), vol. \textbf{23}, Oxford Univ. Press, New York, 2000.

\bibitem{Ey} P. Eymard, L’alg\`ebre de Fourier d’un groupe localement compact, Bull. Soc. Math.
France 92 (1964), 181-236.
%
%

%
%


\bibitem{Fei} H. G. Feichtinger, \textit{On a New Segal Algebra}, Monatsh. Math. 
\textbf{92}(1981), 269-289.


\bibitem{Fo} G.B. Folland, \textit{A Course in Abstract Harmonic Analysis,} 
Studies in Advanced
Mathematics, CRC Press, Boca Raton, FL, 1995.



\bibitem{Jak} M. S. Jakobsen, \textit{On a (No Longer) New Segal Algebra: A Review of the Feichtinger Algebra}, J. Fourier Anal. Appl. \textbf{24} (2018), no. 6, 1579-1660. 

\bibitem{KA} E. Kaniuth, A.T. Lau and A. Ulger, \textit{The Rajchman algebra $B_0(G)$ of a locally compact group $G,$}  Bull. Sci. Math. \textbf{140} (2016), no. 3, 273-302. 

\bibitem{Lamona} R. Lakshmi Lavanya, \textit{On the Fourier transform on function algebras on locally compact Abelian groups,} Monatsh. Math. (4) \textbf{184} (2017), 597-610.

%
%

\bibitem{ReSt} H. Reiter and J. D. Stegeman, \textit{Classical Harmonic
 Analysis and Locally Compact Groups,} LMS Monographs, New series 
 \textbf{22}, Clarendon Press, Oxford, 2000.
 
 \bibitem{Sp}  N. Spronk, \textit{Operator space structure on Feichtinger's Segal algebra,}  J. Funct. Anal. \textbf{248} (2007), no. 1, 152-174.
 
 \bibitem{Wa} M.E. Walter, \textit{Group duality and isomorphisms of Fourier and Fourier-Stieltjes algebras
from a W*-algebra point of view,} Bull. Amer. Math. Soc. \textbf{76} (1970), 1321-1325.

\bibitem{We} J.G. Wendel, \textit{Left centralizers and isomorphisms of group algebras,} Pacific J. Math. \textbf{2} (1952), 251-261.

%

%
%
%
\end{thebibliography}


\end{document}